\documentclass[12pt,a4paper]{article}
\usepackage{amsfonts}
\usepackage{amsbsy}
\usepackage{amssymb}
\usepackage{amsmath}
\usepackage{amsthm}
\usepackage{graphicx}
\usepackage{amsmath}
\newtheorem{theorem}{Theorem}

\newcommand{\titleart}[1]{\bigskip \begin{center}
\large \textbf{#1}\end{center}}
\vspace{-6mm}
\newcommand{\autorart}[1]{\begin{center} \large
\textsf{#1}\end{center} \vspace{-6mm} }
\newcommand{\coord}[1]{\begin{center} \small
\textit{#1}\end{center} \vspace{-6mm} }
\newcommand{\mail}[1]{\begin{center} \small
\textit{#1}\end{center} \normalsize}

\begin{document}
\titleart{On the Ramsey numbers for paths and generalized Jahangir graphs $J_{s,m}$}
\autorart{Kashif Ali, E. T. Baskoro, I. Tomescu}
\coord{COMSATS Institute of Information Technology, Lahore, Pakistan.\\
Combinatorial Mathematics Research Division,
 Institut Teknologi Bandung, Indonesia.\\
Faculty of Mathematics and Computer Sciences,
University of Bucharest,
Str. Academiei, 14,
010014 Bucharest, Romania.}

\mail{akashifali@gmail.com, ebaskoro@math.itb.ac.id, ioan@fmi.unibuc.ro}
\medskip

\begin{abstract}
\noindent For given graphs $G$ and $H,$ the \emph{Ramsey number} $R(G,H)$ is
the least natural number $n$ such that for every graph $F$ of order
$n$ the following condition holds: either $F$ contains $G$ or the
complement of $F$ contains $H.$ In this paper, we determine the
Ramsey number of paths versus  generalized Jahangir graphs. We also
derive the Ramsey number $R(tP_n,H)$, where $H$ is a generalized
Jahangir graph $J_{s,m}$ where $s\geq2$ is even, $m\geq3$ and $t\geq1$
is any integer.
\end{abstract}

\begin{flushleft}
\end{flushleft}

\section{Introduction}

The study of Ramsey numbers for (general) graphs have received
tremendous efforts in the last two decades, see few related papers
\cite{BBS:05}-{\cite{CZZ:04}, \cite{HasB:04,SB:01} and a nice survey
paper \cite{Rad:04}.\\

\noindent Let $G(V,E)$ be a graph with  vertex-set $V(G$) and
edge-set $E(G)$. If $xy \in E(G)$ then $x$ is called {\em adjacent}
to $y$, and $y$ is a {\em neighbor} of $x$ and vice versa. For any
$A \subseteq V(G)$, we use $N_A(x)$ to denote the set of all
neighbors of $x$ in $A$, namely $N_A(x) = \{ y \in A | xy \in
E(G)\}$. Let $P_{n}$ be a path with $n$ vertices, $C_n$ be a cycle
with $n$ vertices,  $W_k$ be a wheel of $k+1$ vertices, i.e., a
graph consisting of a cycle $C_k$ with one additional vertex
adjacent to all vertices of $C_k$.
 For $s,m \geq 2$, the {\em generalized Jahangir graph} $J_{s,m}$ is a
graph on $sm+1$ vertices i.e.,
 a graph consisting of a cycle $C_{sm}$ with one additional vertex which is adjacent
to $m$ vertices of $C_{sm}$ at distance $s$ to each other on
$C_{sm}$.\\

\noindent Recently, Surahmat and Tomescu \cite{ST:06} studied the Ramsey number
of a combination of paths $P_n$ versus  $J_{2,m}$, and obtained the following result.\\

\noindent {\bf Theorem~A}. \cite{ST:06}.\\
 {\it $R(P_n,J_{2,m})=\left\{
\begin{array}{ll}
6&\mbox{if $(n,m)=(4,2),$}\\
n+1& \mbox{if $m=2$ and $n\geq 5,$}\\
 n+m-1& \mbox{~if $m\geq 3$ and $n\geq (4m-1)(m-1)+1.$}
\end{array}\right.$}\\

\noindent
For the Ramsey number of $P_n$ with respect to wheel
$W_m$, Surahmat and Baskoro \cite{BBS:05} showed the following result.\\

\noindent {\bf Theorem B}. \cite {BBS:05}.\\

{\it $R(P_n,W_m)=\left\{
\begin{array}{ll}

2n-1& \mbox{if $m\geq 4$ is even and $n\geq \frac{m}{2}(m-2),$}\\
3n-2& \mbox{~if $m\geq 5$ is odd and $n\geq \frac{m-1}{2}(m-3).$}
\end{array}\right.$}\\

\noindent In this paper, we determine the Ramsey numbers involving
paths $P_n$ and  generalized Jahangir graphs $J_{s,m}.$ We also find
the Ramsey number $R(tP_n,H)$, where $H$ is a generalized Jahangir graph
$J_{s,m}$ where $s\geq2$ is even,
$m\geq3.$ In the following section we prove our main results.\\
\section{Main Results}
\begin{theorem}
 \it For even $s\geq2$ and $m\geq3$,
$R(P_n,J_{s,m})=n+\frac{sm}{2}-1$, where $n\geq
(2sm-1)(\frac{sm}{2}-1)+1$.
\end{theorem}
\noindent {\bf Proof}.\\
 Let $G= K_{n-1}\bigcup
K_{\frac{sm}{2}-1}$. We have $R(P_n,J_{s,m})\geq n+\frac{sm}{2}-1$
since $P_n \not\subseteq G$ and $J_{s,m} \not\subseteq
\overline{G}$. It remains to prove that $R(P_n,J_{s,m})\leq
n+\frac{sm}{2}-1$.  Let $F$ be a graph of order $n+\frac{sm}{2}-1$
and containing no path $P_n$, we will show that
$\overline{F}\supseteq J_{s,m}$. Let $L_1=l_{1,1},l_{1,2}, \ldots,
l_{1,k}$ be the longest path in $F$ and so $k\leq n-1$. If $k=1$ we
have $\overline{F}\cong K_{n+\frac{sm}{2}-1},$ which contains
$J_{s,m}$. Suppose that $k\geq2$ and
$J_{s,m}\not\subseteq\overline{F}.$ We have $zl_{1,1},zl_{1,k}\notin
E(F)$ for each $z\in
V_1=V(F)\backslash V(L_1)$. We distinguish two cases:\\

\noindent{\bf Case 1.} $k\leq 2sm-1$. Let
$L_2=l_{1,2},l_{2,2},\ldots,l_{2,t}$ be a longest path in $F[V_1].$
It is clear that $1\leq t\leq k.$ If $t=1$ then the vertices in
$V_1$ induce a subgraph having only isolated vertices. In this case we
shall add an edge $uv$ to $F$, where $u,v\in V_1$ and denote
$L_2=u,v.$ In this way we can define inductively the system of paths
$L_1,L_2,\ldots,L_{\frac{sm}{2}-1}$ such that $L_i$ is a longest
path in $F[V_{i-1}]$, where $V_{i-1}=V(F)\setminus
\bigcup_{j=1}^{i-1}V(L_j)$ or an edge added to $F$ as above. By
denoting the set of remaining vertices by $B$, we have $|B|\geq
n+\frac{sm}{2}-1-(\frac{sm}{2}-1)(2sm-1)\geq \frac{sm}{2}\geq 3$ since
$s\geq2$ and $m\geq 3$. Let $x,y,z\in B$ be three distinct vertices
which are not in any $L_j$ for $j=1,2,\ldots,\frac{sm}{2}-1$.
Clearly, $x,y,z$ are not adjacent to all endpoints of these $L_j$.
If $F_1$ denotes the graph $F$ or the graph $F$ plus some edges
added in the process of defining the system of paths, it follows
that the endpoints of these $L_j$ induce in $\overline{F_1}$ a
complete graph $K_{sm-2}$ minus a matching having at most
$\frac{sm}{2}-1$ edges if some of the endpoints of same $L_j$ are
adjacent in $F_1$. Since $x,y,z$ are not adjacent to all endpoints
of these $L_j$ it is easy to see that vertices $x,y,z$ and endpoints
of the paths $L_j$ form a
$J_{s,m}\subseteq\overline{F_1}\subseteq\overline{F}$.\\

\noindent{\bf Case 2.} $k > 2sm-1$. In this case we define
$\frac{sm}{2}-1$ quadruple of consecutive vertices of $L_1$ as
follows:

\begin{equation*}
\begin{array}{lll}
C_1 &=& \{l_{1,2},l_{1,3},l_{1,4},l_{1,5}\},\\
C_2 &=& \{l_{1,6},l_{1,7},l_{1,8},l_{1,9}\},\\
& \vdots & \\
C_{\frac{sm}{2}-1} &=&
\{l_{1,2sm-6},l_{1,2sm-5},l_{1,2sm-4},l_{1,2sm-3}\}.
\end{array}
\end{equation*}

\noindent Let $Y=V(F)\setminus V(L_1)$. We have
$|Y|=n+\frac{sm}{2}-1-k\geq\frac{sm}{2}$ since $k\leq n-1$. Hence we
can consider $\frac{sm}{2}$ distinct elements in $Y:
y_1,y_2,\ldots,y_{\frac{sm}{2}}$ and $\frac{sm}{2}-1$ pairs of
elements $Y_i=\{y_i,y_{i+1}\}$ for $i=1,\ldots,\frac{sm}{2}-1$. By
the maximality of $L_1$ it follows that for each
$i=1,\ldots,\frac{sm}{2}-1$ at least one vertex in $C_i$ is not
adjacent to any vertex in $Y_i$. Denote by $c_i$ the vertex in $C_i$
which is not adjacent to any vertex in $Y_i$ for
$i=1,\ldots,\frac{sm}{2}-1$. We have $\overline F \supseteq
J_{s,m}$, where $J_{s,m}$ consists of the cycle $C_{sm}$ having
$V(C_{sm})=\{y_1,c_1,y_2,c_2,\ldots,y_{\frac{sm}{2}-1},
c_{\frac{sm}{2}-1},y_{\frac{sm}{2}},l_{1,k}\}$ and the hub
$l_{1,1}$.\qed \vspace{0.4cm}

\begin{theorem} For odd $s\geq3$,\\
\it $R(P_n,J_{s,m})=\left\{
\begin{array}{ll}
2n-1&\mbox{if $n\geq\frac{sm}{2}(sm-2)$, and $m\geq2$ is even,}\\
2n& \mbox{if $n\geq\frac{sm-1}{2}(sm-1)$, and $m\geq3$ is odd.}
\end{array}\right.$\\
\end{theorem}

\noindent{\bf Proof.}\\

\noindent To show the lower bound, consider graphs $2K_{n-1}$ and
$K_1\cup2K_{n-1}$ for the first and
second cases of Theorem  respectively.\\

\noindent For the  reverse inequality, firstly we will prove the
result for the first case of Theorem. Let $F$ be a graph of order
$2n-1$ containing no path $P_n$ where $n\geq\frac{sm}{2}(sm-2)$. We
will show that $\overline{F}\supseteq J_{s,m}$. Since $F$ does not
contain $P_n$, by Theorem B, $\overline{F}$ will contain a wheel
$W_{sm}$, and  so clearly $ \overline{F}\supseteq J_{s,m}$.\\

\noindent For the second case, \noindent to prove
$R(P_n,J_{s,m})\leq 2n $  let $F$ be a graph on $2n$
vertices containing no $P_n$. Let $L_1 = (
l_{11},l_{12},\cdots,l_{1k-1},l_{1k})$ be a longest path in $F$ and
so $k\leq n-1$. If $k=1$ we have $\overline{F}\simeq K_{2n}$, which
contains $J_{s,m}$. Suppose that $k\geq2$ and $\overline{F}$ does
not contain $J_{s,m}$. Obviously, $ zl_{11},zl_{1k}$ are not in
$E(F)$ for each $z\in V_1$, where $V_1=V(F)\setminus V(L_1)$. Let
$L_2 = ( l_{21},l_{22},\cdots,l_{2t-1},l_{2t})$ be a longest path in
$F[V_1]$. It is clear that $1\leq t\leq k$. Let $V_2=V(F)\setminus
(V(L_1)\cup V(L_2))$. We distinguish three cases.\\

\noindent {\bf Case 1} : $k < {sm-1}$.  If $t=1$ then the vertices
in $V_1$ induce a subgraph having only isolated vertices. In this
case we shall add an edge $uv$ to $F$, where $u,v\in V_1$ and denote
$L_2=u,v$. In this way we can define inductively the system of paths
$L_1,L_2,\cdots, L_\frac{sm-1}{2}$ such that $L_i$ is a longest path
in $F[V_{i-1}]$, where $V_{i-1}=V(F)\setminus \bigcup _{j=1}^{i-1}
V(L_j)$ or an edge added to $F$ as above. If $F_1$ denotes the graph
$F$ or the graph $F$ plus some edges added in the process of
defining the system of paths, it follows that endpoints of these
$L_j$ , where $j=1,2,\cdots,\frac{sm-1}{2}$ induce in
$\overline{F_1}$ a complete graph $K_{sm-1}$ minus a matching having
at most $\frac{sm-1}{2}$ edges if some of the endpoints of same
$L_j$ are adjacent in $F_1$. Since for $s,m\geq3$ there exist at
least two vertices $x,y$ which are not adjacent to all endpoints of
these $L_j$.  Thus, it is easy to see that vertices $x,y$ together
with all endpoints of paths $L_j$  form a  $J_{s,m}\subseteq\overline{F_1}\subseteq\overline{F}$.\\

\noindent {\bf Case 2:} $k\geq sm-1$ and $t\geq sm-1$. For
$i=1,2,\cdots,\frac{sm-3}{2}$ define the couples $A_i$ in path $L_1$
as follows:\\

{\it $A_i=\left\{
\begin{array}{ll}
\{l_{1i+1},l_{1i+2}\} &\mbox {for  i  odd,}\\
\{l_{1k-i},l_{1k-i+1}\} & \mbox{for i  even. }\\

\end{array}\right.$}\\

\noindent Similarly, define couples $B_i$ in path $L_2$ as follows:

{\it $B_i=\left\{
\begin{array}{ll}
\{l_{2i+1},l_{2i+2}\} &\mbox {for  i  odd,}\\
\{l_{2t-i},l_{2t-i+1}\} & \mbox{for i  even. }\\

\end{array}\right.$}\\

\noindent Since $t\leq k\leq n-1$ and $|F| =2n$, there exist at
least two vertices $x,y$ which are not in $L_1\cup L_2$.
 Since $L_1$ is a longest path in $F$, there exists one vertex of
$A_i$ for each $i$, say $a_i$ which is not adjacent with $x$.
 Similarly, since $L_2$ is a longest path in $V(F)\setminus V(L_1)$
there must be one vertex, say $b_i$, in couple $B_i$ which is not
adjacent to $x$ for each $i$.  By maximality of path $L_1$, $b_ia_i$
 and $a_ib_{i+1}$ are not in $E(F)$ for each $i$. Thus
$\{l_{11},b_1,a_1,b_2,a_2,\cdots,b_\frac{sm-3}{2},a_\frac{sm-3}{2},l_{2t},y\}$
will form a cycle $C_{sm}$ in $\overline{F}$  and since $x$ is
adjacent with at least $sm-1$ vertices of cycle $C_{sm}$ in
$\overline{F}$, we have a
subgraph in $\overline{F}$ which contain $J_{s,m}$ , so  $J_{s,m}\subseteq\overline{F}$.\\

\noindent{\bf Case 3:} $k\geq sm-1$ and $t<sm-1$. Since $k\leq n-1$
( $F$ has no $P_n$), $V_1$ will have at least $n+1$ vertices. Then,
we can define the same process as in Case 1, since
$n+1-(sm-2)\frac{sm-1}{2}\geq\frac{sm+1}{2}\geq5$.
\qed \vspace{0.4cm}

\noindent In the following theorem we derive Ramsey number $R(tP_n,J_{s,m})$ for
any integer $t\geq1,$ even $s$ and $m\geq3,$ where $n$ is large enough
with respect to $s$ and $m$ as follows.

\begin{theorem}
$R(tP_n,J_{s,m})= tn+\frac{sm}{2}-1$ if $n\geq(\frac{sm}{2}-1)(2sm-1)+1$,\\
 $s\geq2$ is even,
 $m\geq3$ and $t$ is any positive integer.

\end{theorem}
\noindent { \bf Proof of Theorem 3.}
\noindent  Since graph $G = K_{\frac{sm}{2}-1}\cup K_{tn-1}$
contains no $tP_n$ and $\overline{G}$ contains no $J_{s,m}$, then
$R(tP_n,J_{s,m})\geq tn+\frac{sm}{2}-1$. For proving the upper
bound, let $F$ be a graph of order $tn+\frac{sm}{2}-1$ such that
$\overline{F}$ contains no  $J_{s,m}$.  We will show that $F$
contains $t P_n$. We use  induction on $t$. For $t=1$ this is true
from Theorem 1. Now, let assume that the theorem is true for all
$t^{'} \leq t-1$. Take any graph $F$ of $tn+\frac{sm}{2}-1$ vertices
such that its complement contains no $J_{s,m}$. By the induction
hypothesis, $F$ must contain $t-1$ disjoint copies of $P_n$. Remove
these copies from $F$, then by Theorem 1 the  subgraph $F[H]$ on
remaining vertices will induce another $P_n$ in $F$ since
$\overline{F} \not\supseteq J_{s,m}$, so
$\overline{F[H]}\not\supseteq J_{s,m}$. Therefore $F \supseteq
tP_n$. The proof is complete. \qed

\end{document}